\documentclass[english,12pt]{article}
\usepackage{graphicx}
\usepackage{pifont,amsmath,amsthm,geometry,txfonts,fleqn,url,color,soul}
\usepackage{amsbsy,amscd,amssymb,graphicx,epsfig,babel}
\usepackage{authblk}
\usepackage{times}
\usepackage[colorlinks=true]{hyperref}
\usepackage{amssymb}

\newtheorem{proposition}{Proposition}[section]
\newtheorem{corollary}[proposition]{Corollary}
\newtheorem{lemma}[proposition]{Lemma}
\newtheorem{theorem}[proposition]{Theorem}

\newcommand\s[1]{\{#1\}}

\newcommand{\keywords}{\textbf{Keywords:~}}

\begin{document}


\title{Hypertrees and their host trees: a survey}

\author{Pablo De Caria Di Fonzo\footnote{Webpage:\href{www.mate.unlp.edu.ar/~pdecaria}{www.mate.unlp.edu.ar/~pdecaria}}}

%


\affil{pdecaria@mate.unlp.edu.ar}
\affil{CONICET/ Centro de Matem\'{a}tica de La Plata, Universidad Nacional de La Plata, La Plata, Argentina.}

\date{}

\maketitle

\begin{abstract}
A hypergraph $\mathcal{H}=(V,\mathcal{E})$ is a hypertree if it admits a tree $T$ with vertex set $V$ such that every edge of $\mathcal{H}$ induces a subtree of $T$. A tree like that is called a host tree. Several characterizations and properties of hypertrees have been discovered over the years. However, the interest in the structure of their host trees was weaker and restricted to particular scenarios where they arise, like the clique tree of chordal graphs. In that special case, the proofs of most characteristics of clique trees that exist in the literature rely significantly on the structural properties of chordal graphs. The purpose of this work is the study of the properties of the host trees of hypertrees in a more general context and have them described in a single place, giving simpler proofs for known facts, generalizing others and introducing some new concepts that the author considers that are relevant for the study of the topic. Particularly, we will determine what edges can be found in some host tree of a hypertree, and how these edges must be combined to form a host tree, with an emphasis in tools like the basis and the completion of a hypergraph, and the concept of equivalent hypergraphs.
\end{abstract}

\keywords{Hypergraph, hypertree, host tree, chordal graph.}

%

\section{Introduction}

A \emph{graph} $G$ is a pair $(V,E)$, where $V$ is a set whose elements are called the \emph{vertices} and $E=\s{e_i}_{i\in I}$, where $I$ is a finite set and, for every $i\in I$, $e_i$ is a subset of $V$ with one or two elements, which is called an \emph{edge}. In this work, we will only deal with \emph{simple graphs}, which are those that do not have edges with just one vertex (\emph{loops}) or multiple edges (edges that have exactly the same vertices). Graphs are one of the major objects of study in Discrete Mathematics, but the limitation that edges correspond to sets with one or two elements gave rise to the more general concept of hypergraph.

A \emph{hypergraph} $\mathcal{H}$ is also defined as a pair $(V,\mathcal{E})$, where $V$ is again the vertex set, $\mathcal{E}=\s{e_i}_{i\in I}$, $I$ is a finite set and, for every $i\in I$, $e_i$ is a nonempty subset of $V$, which is called an edge (or \emph{hyperedge}). When no subset is repeated, case in which we say that the hypergraph is \emph{simple}, $\mathcal{E}$ is just a subset of the power set of $V$. As for notation, the nature of edges as subsets of the vertex set will sometimes make it convenient to refer to them using set notation (upper case letters) even if the usual thing is to use lower case letters to denote an edge. Both ways are used in this work; the context will determine our preference. We may also use the notation $V(\mathcal{H})$ and $E(\mathcal{H})$ to refer to the vertex set and family of edges of $\mathcal{H}$, respectively. Hypergraphs have several applications, like for data analysis, system modeling \cite{modeling}, image retrieval \cite{retrieval} and bioinformatics \cite{bioinf}, to cite some examples.

\begin{figure}[h]\label{hypex}
\begin{center}
\scalebox{0.75}{\includegraphics{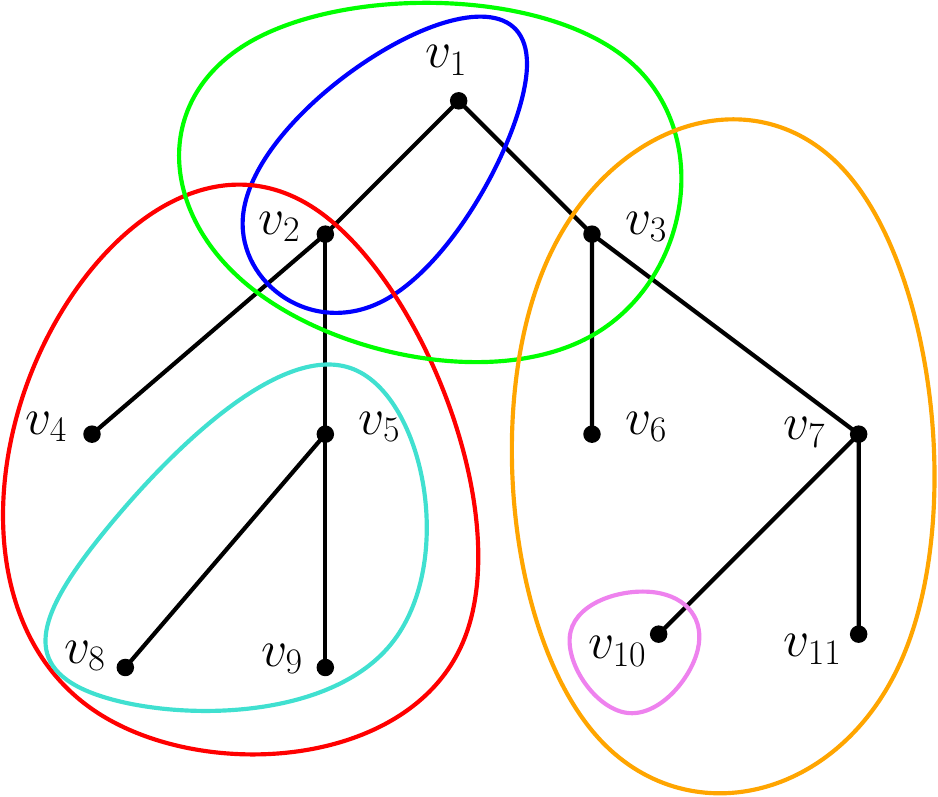}}
\caption{A hypertree with eleven vertices and six edges, that are represented as closed curves. The tree that appears in the figure is a host tree of it.}
\end{center}
\end{figure}

In this paper, we are interested in a special type of hypergraph. A hypergraph $\mathcal{H}$ is a \emph{hypertree} if there exists a tree $T$ whose vertex set is $V(\mathcal{H})$ and, for every edge $e$ of $\mathcal{H}$, the subgraph of $T$ induced by $e$ is a subtree. Such a tree will be called a \emph{host tree} of $\mathcal{H}$ (see Figure \ref{hypex} for an example).

This is not the only way a hypertree can be described. We say that a simple graph $G$ is \emph{chordal} if it does not have induced cycles of length four or greater. Given a hypergraph $\mathcal{H}$, the \emph{line graph} of $\mathcal{H}$, or $L(\mathcal{H})$, has all the edges of $\mathcal{H}$ as vertices in such a way that two different edges $e_1$ and $e_2$ of $\mathcal{H}$ are adjacent in $L(\mathcal{H})$ if and only if their intersection is not empty. We also say that $\mathcal{H}$ is \emph{Helly} if for every instance of edges $e_1$, $e_2$,..., $e_k$ of $\mathcal{H}$ that are pairwise intersecting, the intersection of all of them is not empty.

\begin{theorem} \label{arbor}\cite{arbor}
A hypergraph $\mathcal{H}$ is a hypertree if and only $\mathcal{H}$ is Helly and $L(\mathcal{H})$ is chordal.
\end{theorem}

Given a hypergraph $\mathcal{H}$, the \emph{dual hypergraph} $D\mathcal{H}$ is the hypergraph whose vertices are the edges of $\mathcal{H}$ and that, for every vertex $v$ of $\mathcal{H}$, has the edge $D_v:=\s{e\in \mathcal{E}:~v\in e}$. It is simple to verify that the hypergraph $DD\mathcal{H}$ is isomorphic to $\mathcal{H}$, making every vertex $v$ of $\mathcal{H}$ correspond to the vertex $D_v$ of $DD\mathcal{H}$ and every edge $F$ of $\mathcal{H}$ correspond to the edge $D_F$ of $DD\mathcal{H}$. A hypergraph is said to be \emph{conformal} when its dual is Helly.

Define a \emph{dual hypertree} as a hypergraph $\mathcal{H}$ for which there exists a tree $T$ whose vertices are the edges of $\mathcal{H}$ and such that, for every vertex $v$ of $\mathcal{H}$, the edges of $\mathcal{H}$ that contain $v$ induce a subtree. It is straightforward that a hypergraph is a dual hypertree if and only if its dual hypergraph is a hypertree. Dual hypertrees play an important role in the theory of desirable properties of relational database schemes \cite{algor,fagin}.

A \emph{hypercycle} of a hypergraph $\mathcal{H}$ is a cycle of $L(\mathcal{H})$. By the characterization of hypertrees we have seen, if $e_1e_2...e_ke_1$ is a hypercycle of length greater than or equal to four of a hypertree $\mathcal{H}$, the chordality of $L(\mathcal{H})$ implies that there are two edges of $\mathcal{H}$ that are in the hypercycle, are not consecutive in it and their intersection is not empty. This acyclicity that hypertrees satisfy can be redefined to characterize dual hypertrees.

Given a hypercycle $C:e_1e_2...e_ke_1$ of the hypergraph $\mathcal{H}$, a \emph{chord} of $C$ is an edge $e$ of $\mathcal{H}$ such that $e_i\cap e_{i+1}\subseteq e$ (mod $k$) for at least three values of $i$ between 1 and $k$. In this context, dual hypertrees can be characterized in the following way:

\begin{theorem} \cite{fagin2,fagin3}
A hypergraph is a dual hypertree if and only if it is conformal and every hypercycle of length greater than or equal to three has a chord.
\end{theorem}

Finally, there is one more characterization:

\begin{theorem} \cite{fagin2}
A hypergraph $\mathcal{H}=(V,\mathcal{E})$ is a dual hypertree if and only if $\mathcal{E}$ can be emptied through the application of the following operations:

\begin{enumerate}
  \item Removing an edge that is contained in another.
  \item If $v$ is a vertex that appears only in one edge $e$, remove $v$ from $e$.
\end{enumerate}

\end{theorem}

Whereas host trees have been important to implement the applications of hypertrees in general, there has not been any work that mainly focuses on them regardless of the uses they can have. They have received the most attention in special cases, like in the context of chordal and dually chordal graphs. A graph $G$ is chordal if and only if it has a \emph{clique tree}, that is, a tree $T$ whose vertices are the maximal cliques of $G$ and such that, for every $v\in V(G)$, the set $\mathcal{C}_v$ of maximal cliques of $G$ that contain $v$ induces a subtree of $T$ \cite{gavril}. In other words, clique trees are the host trees of the dual of the hypergraph of maximal cliques, so chordal graphs can be characterized as the graphs whose hypergraph of maximal cliques is a dual hypertree.

Denote the set of maximal cliques of a graph $G$ by $\mathcal{C}(G)$. The \emph{clique graph} of a graph $G$, or $K(G)$, is the line graph of $\mathcal{C}(G)$. A graph is said to be \emph{dually chordal} if it is the clique graph of a chordal graph. A necessary and sufficient condition for a graph $G$ to be dually chordal is the existence of a tree $T$ that has the same vertices as $G$ and such that, for every maximal clique $C$ of $G$, $C$ induces a subtree of $T$ \cite{brands}. A tree like that is usually called a \emph{compatible tree} and, in connection to hypertrees, it is a host tree of the hypergraph of maximal cliques of $G$. Thus, a graph is dually chordal if and only if its maximal cliques make a hypertree. Compatible trees can also be characterized as those trees with the same vertices as the graph and such that the closed neighborhood of every vertex induces a subtree. Consequently, a graph $G$ is dually choral if and only if the hypergraph of the closed neighborhoods of the vertices of $G$ is a hypertree.

For clique trees and compatible trees, many properties of them are known, which were proved using the structural properties of chordal \cite{reduced,enum} and dually chordal graphs. The purpose of this paper is to conduct a general study of the host trees of hypertrees and their structure, providing a survey of their properties and giving simpler and elementary proofs, that will have the facts about clique trees and compatible trees as corollaries.

In Section 2, we introduce the notion of equivalent hypergraphs, that is, hypergraphs that have the same host trees. Particularly, we define the major two equivalent hypergraphs that any hypergraph has, namely, the completion and the basis. We also use the notion of equivalence to define basic hypertrees, which are connected to basic chordal graphs.

In Section 3, we take advantage of all what was proved in the previous section to find conditions for an edge to be in the host tree of a hypergraph and use it to describe the structure of host trees, also exploring how this translates to the tree representations of chordal and dually chordal graphs. In particular, we give a new characterization of hypertrees.

In Section 4, we briefly discuss how host trees can be characterized as the maximum weight spanning trees of a weighted graph, although not every set of maximum weight spanning trees of a weighted graph are exactly the host trees of some hypergraph.

In Section 5, we present the conclusions.

\section{Equivalent hypertrees}

In this section, we explore the notion of equivalent hypertrees as an opportunity to introduce new special hypertrees that will facilitate the study of the structure of host trees.

We say that two hypergraphs $\mathcal{H}$ and $\mathcal{H}'$ on the same vertex set $V$ are \emph{equivalent} if they have the same host trees. As a consequence of this definition, if $\mathcal{H}$ is not a hypertree, then neither is $\mathcal{H}'$ for them to be equivalent.

A union $\cup_{i=1}^n{F_i}$ is said to be \emph{connected} if the line graph induced by these sets is connected. The symbol $\biguplus$ will be used to emphasize that a union is connected.

Here we list some basic operations that can be applied to a hypergraph to obtain another equivalent hypergraph. Their justifications are quite trivial, but we include details for completion and to simplify ensuing proofs.

\begin{proposition} \label{operations}
  Let $\mathcal{H}$ be a hypergraph. The following operations applied to $\mathcal{H}$ yield an equivalent hypergraph.

  \begin{enumerate}
    \item Removing/adding to $\mathcal{H}$ an edge with just one element.
    \item Removing/adding to $\mathcal{H}$ an edge with all the vertices of $\mathcal{H}$.
    \item Adding to $\mathcal{H}$ an edge with the same vertices as an edge already in $\mathcal{H}$.
    \item Removing an edge of $\mathcal{H}$ for which there exists another edge with the same vertices.
    \item Adding to $\mathcal{H}$ an edge that is the nonempty intersection of some edges of $\mathcal{H}$.
    \item Adding to $\mathcal{H}$ an edge that is the union of two edges of $\mathcal{H}$ that are not disjoint.
    \item Adding to $\mathcal{H}$ an edge that is the connected union of some edges of $\mathcal{H}$.
    \item Removing from $\mathcal{H}$ an edge that is the intersection of other edges of $\mathcal{H}$.
    \item Removing from $\mathcal{H}$ an edge that is the connected union of other edges of $\mathcal{H}$.
    \item Any combination of the previous operations.
  \end{enumerate}
\end{proposition}

\begin{proof}
  1 and 2 are true because those edges induce a subtree of any tree with vertex set $V(\mathcal{H})$, so they do not create additional restrictions for a host tree. 3 and 4 are true because those operations do not change the subsets of $V(\mathcal{H})$ that must induce a subtree for a tree $T$ to be a host tree.

  Let $A_1$ and $A_2$ be edges of $\mathcal{H}$ that are not disjoint and let $T$ be a host tree of $\mathcal{H}$. Let $B=A_1\cup A_2$. If $B=A_1$ or $B=A_2$, then it is clear that $T[B]$ is a subtree. Otherwise, let $x$ and $y$ be any two different elements of $B$. Consider an element $z$ in the intersection of $A_1$ and $A_2$. Then, $T[x,z]$ and $T[z,y]$ are paths in $T[A_1]$ or in $T[A_2]$. Thus, the union of $T[x,z]$ and $T[z,y]$ gives a walk from $x$ to $y$ in $T[B]$. We conclude that $T[B]$ is connected and hence a subtree. This argument implies that every host tree of $\mathcal{H}$ is a host tree of the new hypergraph, and it is trivial that the converse is also true. Therefore, the hypergraphs are equivalent, which proves 6.

  7 follows from 6 noting that every connected union can be obtained from successive unions of two sets with nonempty intersections.

  Suppose now that $A_1,...,A_k$ are edges of $\mathcal{H}$ with nonempty intersection and $T$ is a host tree of $\mathcal{H}$. Let $B=\cap_{i=1}^k{A_i}$. If $B$ has a single element, then $T[B]$ is a subtree. Otherwise, let $x$ and $y$ be any two different elements of $B$. Since, for $1\leq i\leq k$, $T[A_i]$ is a subtree of $T$, $T[x,y]$ is a path in $T[A_i]$. Thus, $T[x,y]$ is a path of $T[B]$. We conclude that $T[B]$ is connected and hence a subtree. This argument implies that every host tree of $\mathcal{H}$ is a host tree of the new hypergraph, and it is trivial that the converse is also true. Therefore, the hypergraphs are equivalent, which proves 5.

  If we apply 8 or 9 to $\mathcal{H}$ to get a new hypergraph, then we can next apply 5 or 7 (respectively) to get back to $\mathcal{H}$, from which we can infer that 8 and 9 indeed yield equivalent hypergraphs.

  10 is straightforward since each operation applied yields a hypergraph equivalent to the previous one and, by transitivity, to $\mathcal{H}$.
\end{proof}



We will denote by $Simp(\mathcal{H})$ the hypergraph obtained from $\mathcal{H}$ by deleting repeated edges. By 4 and 10 of Proposition \ref{operations}, $\mathcal{H}$ and $Simp(\mathcal{H})$ are equivalent.



As for intersections and connected unions, we will use them to define the next equivalent hypergraph. Let the \emph{completion} of the hypergraph $\mathcal{H}$, or $Comp(\mathcal{H})$, be the simple hypergraph with the same vertex set as $\mathcal{H}$ and whose edges are $V(\mathcal{H})$, the unit subsets of $V(\mathcal{H})$ and the proper subsets of $V(\mathcal{H})$ that can be obtained from the edges of $\mathcal{H}$ through the application (in any order and amount) of the operations of intersection and connected union.

By 1, 2, 5, 7 and 10 of Proposition \ref{operations}, $\mathcal{H}$ and $Comp(\mathcal{H})$ are equivalent hypergraphs.

Suppose now that $\mathcal{H}$ is a hypergraph whose edge set is closed under connected unions. Following what is done in \cite{correspondence}, we define the \emph{basis} of $\mathcal{H}$, or $\mathcal{B}(\mathcal{H})$, as the simple hypergraph with the same vertex set as $\mathcal{H}$ and whose edges are the edges of $\mathcal{H}$ that have more than one vertex and cannot be expressed as the connected union of smaller edges of $\mathcal{H}$. As a consequence of this definition, the sets in $E(Simp(\mathcal{H}))$ that have more than one vertex are just those that can be expressed as a connected union of edges of $\mathcal{B}(\mathcal{H})$, and $\mathcal{B}(\mathcal{H})$ is minimal with respect to this property. By 1, 9 and 10 of Proposition \ref{operations}, $\mathcal{H}$ is equivalent to $\mathcal{B}(\mathcal{H})$. Also note that every edge of $\mathcal{H}$ that has two vertices appears in the basis of $\mathcal{H}$, since it cannot be expressed as the connected union of singletons.

If $\mathcal{H}$ does not have its edges closed under connected unions, $Comp(\mathcal{H})$ does, so we will often refer to the edges of $\mathcal{B}(Comp(\mathcal{\mathcal{H}}))$ as \emph{basic sets} of $\mathcal{H}$.

If $\mathcal{H}$ has its edges closed under connected unions, $\mathcal{H}$ and $Comp(\mathcal{H})$ may have a different basis. Consider for example the simple hypergraph $\mathcal{H}$ with vertex set $\s{1,2,3,4}$ and whose edges are the sets $\s{1,2,3}$, $\s{2,3,4}$ and $\s{1,2,3,4}$. The basis of $\mathcal{H}$ has just the sets $\s{1,2,3}$ and $\s{2,3,4}$, while the basis of $Comp(\mathcal{H})$ also has the edge $\s{2,3}$, the intersection of the other two sets. Even in this case, we will use the term basic set to refer to an edge of $\mathcal{B}(Comp(\mathcal{H}))$, regardless of the fact that $\mathcal{H}$ has its own basis.

In any case, every edge of $\mathcal{H}$ is also an edge of $Comp(\mathcal{H})$ and can be expressed as a connected union of basic sets.

For highlighting purposes, we state the following result,




\begin{theorem}
  Let $\mathcal{H}$ be a hypergraph. Then, $\mathcal{H}$ is a hypertree if and only if its completion is a hypertree. Moreover, $\mathcal{H}$ is a hypertree if and only if $\mathcal{B}(Comp(\mathcal{H}))$ is a hypertree.
\end{theorem}

\begin{proof}
  We know that $\mathcal{H}$ is equivalent to $Comp(\mathcal{H})$, from which the first part follows. Additionally, $Comp(\mathcal{H})$ is equivalent to its basis, so it is equivalent to $\mathcal{H}$, and the second part follows.
\end{proof}

The ideas of defining a basis were first applied in \cite{enum}, where the reduced clique hypergraph of a chordal graph is introduced. However, the idea is only applied in the context of chordal graphs, relying in characterizations of them like the existence of a perfect elimination orderings. Here, we will follow a more general approach that works only with the edges of hypergraphs without resorting to special characterizations.

Now, we will consider the case in which $\mathcal{H}$ is a hypertree, to find how the basic sets of $\mathcal{H}$ are related to host trees.

Given a hypergraph $\mathcal{H}$ and a subset $V'$ of $V(\mathcal{H})$, define $I_\mathcal{H}(V')$ as the intersection of all the edges of $\mathcal{H}$ that contain $V'$, or as the whole vertex set of $\mathcal{H}$ if no edge of $\mathcal{H}$ contains $V'$. When $V'$ has just two elements $u$ and $v$, we allow the notation $I_\mathcal{H}(uv)$ for this case. One simple and basic result about this type of sets is the following one:

\begin{proposition} \label{every}
  Let $V_1$ and $V_2$ be sets of vertices of a hypergraph $\mathcal{H}$. Then, $I_{\mathcal{H}}(V_1)\subseteq I_{\mathcal{H}}(V_2)$ if and only if every edge of $\mathcal{H}$ that contains $V_2$ also contains $V_1$.
\end{proposition}

\begin{proof}
  Suppose that every edge of $\mathcal{H}$ that contains $V_2$ also contains $V_1$. If no edge of $\mathcal{H}$ contains $V_2$, then $I_{\mathcal{H}}(V_2)=V(\mathcal{H})$ and the inclusion $I_{\mathcal{H}}(V_1)\subseteq I_{\mathcal{H}}(V_2)$ trivially follows. If there is an edge of $\mathcal{H}$ that contains $V_2$, then both $I_{\mathcal{H}}(V_1)$ and $I_{\mathcal{H}}(V_2)$ are equal to intersections of sets, the first one having all the sets (and possibly more) of the second. Thus, the inclusion $I_{\mathcal{H}}(V_1)\subseteq I_{\mathcal{H}}(V_2)$ follows easily too.

  Suppose now that $I_{\mathcal{H}}(V_1)\subseteq I_{\mathcal{H}}(V_2)$. It follows from the definition that $V_1\subseteq I_{\mathcal{H}}(V_1)$, so $V_1\subseteq I_{\mathcal{H}}(V_2)$ as well. Therefore, by the definition of $I_{\mathcal{H}}(V_2)$, every edge of $\mathcal{H}$ that contains $V_2$ should also contain $V_1$.
\end{proof}

Now we will see that these sets $I_{\mathcal{H}}(V')$ play an important role in relation to host trees of hypertrees. Particularly, we will use them to find new characterizations of the edges of $Comp(\mathcal{H})$ and its basis.

\begin{lemma} \label{-+}
  Let $\mathcal{H}$ be a hypertree, $T$ be a host tree of $\mathcal{H}$, $uv$ be an edge of $T$ and $x$, $y$ be two vertices such that $T[x,y]$ contains $u$ and $v$ and every edge of $\mathcal{H}$ that contains $\s{u,v}$ also contains $\s{x,y}$. Then, $T-uv+xy$ is a host tree of $\mathcal{H}$.
\end{lemma}

\begin{proof}
  Denote the tree $T-uv+xy$ by $T'$ for short.

Every edge of $\mathcal{H}$ that does not contain $\s{u,v}$ induces the same subtree in $T$ and $T'$. On the other side, by the definition of host tree and the hypothesis, the edges of $\mathcal{H}$ that contain $\s{u,v}$ are just the edges of $\mathcal{H}$ that contain $\s{x,y}$, so every edge $F$ of $\mathcal{H}$ that contains $\s{u,v}$ induces the subtree $T[F]-uv+xy$ of $T'$.

Therefore, $T'$ is a host tree of $\mathcal{H}$.

\end{proof}

\begin{proposition} \label{completion}
  Let $\mathcal{H}$ be a hypertree, $T$ be a host tree of $\mathcal{H}$ and $F$ be a subset of $V(\mathcal{H})$ that induces a subtree of every host tree of $\mathcal{H}$ and is not a unit set. Then

  \begin{enumerate}
    \item If $uv$ is an edge of $T$ such that $\s{u,v}\subseteq F$, then $I_\mathcal{H}(uv)\subseteq F$.
    \item $F=\biguplus_{uv\in E(T[F])}{I_\mathcal{H}(uv)}$
     \end{enumerate}

\end{proposition}

\begin{proof}
  Let $uv$ be an edge of the subtree $T[F]$ and let $w$ be any element in $I_\mathcal{H}(uv)$. Suppose without loss of generality that $w$ is in the connected component of $T-uv$ that has $u$, so $u\in T[v,w]$. Since every edge that contains $\s{u,v}$ also contains $w$, we can apply Lemma \ref{-+} to conclude that the tree $T'$ defined by $T'=T-uv+vw$ is a host tree of $\mathcal{H}$. Thus, $T'[F]$ is a subtree. As $F$ contains $u$ and $v$, it contains every vertex of the $uv$-path in $T'$, with $w$ being one of them. Therefore, $w\in F$. It follows that $I_\mathcal{H}(uv)\subseteq F$, so 1 is true. Property 2 is now a direct consequence of 1.
\end{proof}

Proposition \ref{completion} applies particularly to the edges of $Comp(\mathcal{H})$, which must induce a subtree of every host tree of $\mathcal{H}$ due to the equivalence between the two hypergraphs.

\begin{theorem} \label{wabasic}
  Let $\mathcal{H}$ be a hypertree, $T$ be a host tree of $\mathcal{H}$ and $uv$ be an edge of $T$. Then, $I_\mathcal{H}(uv)$ is a basic set of $\mathcal{H}$. What is more, every basic set of $\mathcal{H}$ is of this form.
\end{theorem}

\begin{proof}
  Suppose first that no edge of $\mathcal{H}$ contains $\s{u,v}$, so $I_\mathcal{H}(uv)=V(\mathcal{H})$. Let $A$ and $B$ be the connected components of $T-uv$. Then, every edge of $\mathcal{H}$ different from $V(\mathcal{H})$ is either contained in $A$ or contained in $B$. Every nonempty intersection or connected union of these edges will still be contained in $A$ or $B$, so it also holds that every edge of $Comp(\mathcal{H})$ is contained in $A$ or $B$, which prevents $V(\mathcal{H})$ from being the connected union of some of these sets. Therefore, $V(\mathcal{H})$  is a basic set under the conditions we assumed.

  Suppose now that there exists an edge of $\mathcal{H}$ that contains $\s{u,v}$. If $V(\mathcal{H})$ is the only one satisfying this condition, then we can reason similarly to the previous paragraph considering the other edges of $\mathcal{H}$.

  Consider now the case that $I_\mathcal{H}(uv)$ is strictly contained in $V(\mathcal{H})$. Let $I_\mathcal{H}(uv)=\biguplus_{i=1}^k{F_i}$, a connected union of edges of $Comp(\mathcal{H})$. Since $T[I_\mathcal{H}(uv)]$ is a subtree that contains the edge $uv$, for the union to be connected there has to be some $F_{i}$ that contains $\s{u,v}$. We have by part 1 of  Proposition \ref{completion} that $I_\mathcal{H}(uv)\subseteq F_i$, thus leading to the conclusion that $I_\mathcal{H}(uv)$ and $F_i$ are equal. Therefore, $I_\mathcal{H}(uv)$ satisfies the condition to be a basic set of $\mathcal{H}$.

  Conversely, suppose that $B$ is a basic set of $\mathcal{H}$. Then, $T[B]$ is a subtree. By part 2 of Proposition \ref{completion},  we have that $B=\biguplus_{xy\in E(T[B])}{I_\mathcal{H}(xy)}$. As a basic set, $B$ cannot be the connected union of smaller sets of  $Comp(\mathcal{H})$. Therefore, $B$ is one of the sets of the union.
\end{proof}

The possibility to obtain the basic sets of a hypertree with the aid of a host tree has algorithmic consequences. Given a hypertree, it is possible to find a host tree for it in linear time \cite{spin}. Then, using the edges of this tree, the sets $I_\mathcal{H}(xy)$ can be found. This is more economic than finding all of $Comp(\mathcal{H})$.

Now we aim at proving that $Comp(\mathcal{H})$ is the largest possible simple hypertree that is equivalent to $\mathcal{H}$.

\begin{theorem}
  Let $\mathcal{H}$ be a hypertree. Then, $E(Comp(\mathcal{H}))$ consists of all the nonempty subsets of $V(\mathcal{H})$ that induce a subtree of every host tree of $\mathcal{H}$. Additionally, if $\mathcal{H}'$ is a hypertree equivalent to $\mathcal{H}$, then $E(Simp(\mathcal{H}'))\subseteq E(Comp(\mathcal{H}))$.
\end{theorem}

\begin{proof}
  Since $Comp(\mathcal{H})$ is equivalent to $\mathcal{H}$, every edge of it induces a subtree of every host tree of $\mathcal{H}$.

  On the other side, if $F$ is a set that is not a unit set and induces a subtree of every host tree of $\mathcal{H}$, then by part 2 of Proposition \ref{completion} and Theorem \ref{wabasic} $F$ is the connected union of basic sets of $\mathcal{H}$. We conclude that $F$, as a connected union of edges of $Comp(\mathcal{H})$ is itself an edge of $Comp(\mathcal{H})$.

 Similarly, if $\mathcal{H}'$ is a hypergraph equivalent to $\mathcal{H}$, then every edge of $\mathcal{H}'$ induces a subtree of every host tree of $\mathcal{H}$ and, by the previous paragraph, is in $E(Comp(\mathcal{H}))$.
\end{proof}

The proof of Proposition \ref{completion} makes it now simpler to demonstrate how the equivalence between hypertrees can be tested using basic sets.

\begin{lemma}\label{hosttcont}
  Let $\mathcal{H}_1$ and $\mathcal{H}_2$ be two hypertrees with the same set $V$ of vertices, such that every host tree of $\mathcal{H}_1$ is also a host tree of $\mathcal{H}_2$. Let $T$ be a host tree of both $\mathcal{H}_1$ and $\mathcal{H}_2$. Then, for every edge $uv$ of $T$, $I_{\mathcal{H}_1}(uv)\subseteq I_{\mathcal{H}_2}(uv)$.
\end{lemma}

\begin{proof}
$I_{\mathcal{H}_2}(uv)$ contains $\s{u,v}$ and induces a subtree of every host tree of $\mathcal{H}_2$, and hence of every host tree of $\mathcal{H}_1$. Apply part 1 of Proposition \ref{completion} to this set, $\mathcal{H}_1$, $T$ and $uv$ to get the desired conclusion.
\end{proof}

\begin{theorem} \label{test}
  Let $\mathcal{H}_1$ and $\mathcal{H}_2$ be two hypertrees with the same set $V$ of vertices. Then, $\mathcal{H}_1$ is equivalent to $\mathcal{H}_2$ if and only if they have the same basic sets. Moreover, if $\mathcal{H}_1$ and $\mathcal{H}_2$ share a common host tree $T$, then they are equivalent if and only if $I_{\mathcal{H}_1}(uv)=I_{\mathcal{H}_2}(uv)$ for every edge $uv$ of $T$.
\end{theorem}

\begin{proof}
  We start proving the second part of the theorem. Suppose that $T$ is a common host tree for $\mathcal{H}_1$ and $\mathcal{H}_2$. If $I_{\mathcal{H}_1}(uv)=I_{\mathcal{H}_2}(uv)$ for every edge $uv$ of $T$, then it follows from Theorem \ref{wabasic} that $Comp(\mathcal{H}_1)$ and $Comp(\mathcal{H}_2)$ have the same basis and hence are equal, so $\mathcal{H}_1$ and $\mathcal{H}_2$ are equivalent.

  Conversely, suppose that $\mathcal{H}_1$ and $\mathcal{H}_2$ are equivalent and let $uv$ be an edge of $T$. Every host tree of $\mathcal{H}_1$ is a host tree of $\mathcal{H}_2$ and every host tree of $\mathcal{H}_2$ is a host tree of $\mathcal{H}_1$, so we can apply Lemma \ref{hosttcont} twice to conclude that $I_{\mathcal{H}_1}(uv) = I_{\mathcal{H}_2}(uv)$.



  We now prove the first part. If $\mathcal{H}_1$ and $\mathcal{H}_2$ have the same basic sets, then we can reason like in the first paragraph of the proof, so $\mathcal{H}_1$ and $\mathcal{H}_2$ are equivalent.

  Conversely, suppose that $\mathcal{H}_1$ and $\mathcal{H}_2$ are equivalent. Let $B$ be a basic set of $\mathcal{H}_1$. We now prove that $B$ is also a basic set of $\mathcal{H}_2$. Let $T$ be a host tree of $\mathcal{H}_1$ and $\mathcal{H}_2$. By Theorem \ref{wabasic}, $B=I_{\mathcal{H}_1}(uv)$ for some edge $uv$ of $T$. Given that we have proved that $I_{\mathcal{H}_1}(uv)=I_{\mathcal{H}_2}(uv)$, we conclude, again by Theorem \ref{wabasic}, that $B$ is a basic set of $\mathcal{H}_2$ as well. Consequently, every basic set of $\mathcal{H}_1$ is also a basic set of $\mathcal{H}_2$. Similarly, every basic set of $\mathcal{H}_2$ is also a basic set of $\mathcal{H}_1$. Therefore, $\mathcal{H}_1$ and $\mathcal{H}_2$ have the same basic sets.


    \end{proof}

 Finally, we now prove that, for a hypertree $\mathcal{H}$, the operations of Proposition \ref{operations} are enough to get all the hypertrees equivalent to it.

 \begin{theorem}
   Let $\mathcal{H}$ and $\mathcal{H}'$ be two hypertrees with the same vertex set. Then, they are equivalent if and only if $\mathcal{H}'$ can be obtained from $\mathcal{H}$ through the operations of Proposition \ref{operations}.
 \end{theorem}

 \begin{proof}
   We know that the operations of Proposition \ref{operations} yield equivalent hypergraphs. For that reason, we only need to prove one of the directions of the theorem.

   Suppose that $\mathcal{H}$ and $\mathcal{H}'$ are equivalent hypertrees. Then, $Comp(\mathcal{H})$ and $Comp(\mathcal{H}')$  have the same basis $\mathcal{B}$. Given that every basic set can be expressed as intersection of edges of $\mathcal{H}$, operation 5 can be applied to $\mathcal{H}$, if necessary,  to get a hypergraph that has all the edges of $\mathcal{B}$. Since every edge of $\mathcal{H}$ is the connected union of edges of $\mathcal{B}$, operation 9 could now be applied to be only left just with $\mathcal{B}$. Next, operation 7 can be applied to have all the edges of $\mathcal{H}'$. Finally, in case of necessity, operation 8 can be applied to remove the basic sets that are not edges of $\mathcal{H}'$, to end up just with $\mathcal{H}'$.
    
 \end{proof}

This proof is only valid for hypertrees because all hypergraphs that share the vertex set and are not hypertrees are technically equivalent without the need of having the same basic sets, so the argument of the proof does not apply. It is true that applying one of the operations to a hypergraph $\mathcal{H}$ that is not a hypertree yields another hypergraph that is not a hypertree, and hence equivalent to $\mathcal{H}$. However, the converse is not true. Consider for example the hypergraphs $\mathcal{H}_1$ and $\mathcal{H}_2$ with vertex set $\s{1,2,3,4,5,6}$ and such that the edges of $\mathcal{H}_1$ are $\s{1,2}$, $\s{2,3}$, $\s{1,3}$, $\s{4,5}$, $\s{4,6}$ and $\s{5,6}$ and the edges of $\mathcal{H}_2$ are $\s{1,2}$, $\s{1,4}$, $\s{2,4}$, $\s{3,5}$, $\s{3,6}$ and $\s{5,6}$. $\mathcal{H}_1$ and $\mathcal{H}_2$ are not hypertrees and their edges are just their basic sets. However, no operation of Proposition \ref{operations} applied to the edges of $\mathcal{H}_1$ allows to get the edge $\s{1,4}$.

It is also interesting to note that, by their definitions, all the operations preserve the basic sets, since they cannot change the completion of the hypergraph.


\vspace{10pt}

  For a hypergraph $\mathcal{H}$, the \emph{2-section} of $\mathcal{H}$, symbolized $2S(\mathcal{H})$, is the graph whose vertex set is $V(\mathcal{H})$, and such that two vertices $u$ and $v$ of $\mathcal{H}$ are adjacent if there exists an edge of $\mathcal{H}$ that contains both vertices. Also define $\mathcal{N}(\mathcal{H})$ as the hypergraph that has the same vertices as $\mathcal{H}$ and that has, for every $v\in V(\mathcal{H})$, the edge $\cup_{F\in D_v}{F}$, the set consisting of the union of all the edges of $\mathcal{H}$ that have $v$. Note that $\mathcal{N}(\mathcal{H})$ is the hypergraph of the closed neighborhoods of the 2-section of $\mathcal{H}$. Define a \emph{basic hypertree} as a hypergraph $\mathcal{H}$ that is a hypertree and such that $\mathcal{H}$ is equivalent to $\mathcal{N}(\mathcal{H})$.

Theorem \ref{test} gives the following necessary and sufficient condition for a hypertree to be basic.

\begin{lemma}
  Let $\mathcal{H}$ be a hypergraph. Then, every host tree of $\mathcal{H}$ is a host tree of $\mathcal{N}(\mathcal{H})$.
\end{lemma}

\begin{proof}
  Let $T$ be a host tree of $\mathcal{H}$ and let $F$ be an edge of $\mathcal{N}(H)$, corresponding to the closed neighborhood of a vertex $x$ in the 2-section of $\mathcal{H}$. Thus, $F$ is the connected union of all the edges of $\mathcal{H}$ that contain $x$, so $F$ is an edge of $Comp(\mathcal{H})$ and $T[F]$ must be a subtree.

  Therefore, $T$ is also a host tree of $\mathcal{N}(\mathcal{H})$.
\end{proof}

\begin{proposition} \label{basicn}
  Let $\mathcal{H}$ be a hypertree and $T$ be a host tree of $\mathcal{H}$ and $\mathcal{N}(\mathcal{H})$. Then, the basic set $I_{\mathcal{N}(\mathcal{H})}(uv)$ consists of all the vertices $x$ satisfying that, for every vertex $y$ such that there exists an edge of $\mathcal{H}$ that contains $\s{u,v,y}$ (alternatively, that contains $I_\mathcal{H}(uv)\cup\s{y}$), there exists an edge of $\mathcal{H}$ that contains $\s{x,y}$.
\end{proposition}

\begin{proof}

If no edge of $\mathcal{H}$ contains $\s{u,v}$, then it follows that no edge of $\mathcal{N}(\mathcal{H})$ contains $\s{u,v}$, either. Thus, $I_{\mathcal{N}(\mathcal{H})}(uv)=V(\mathcal{H})$ and the Proposition is true by default. Suppose from now on that some edge of $\mathcal{H}$ contains $\s{u,v}$.

Let $x\in I_{\mathcal{N}(\mathcal{H})}(uv)$ and let $y$ be a vertex such that there exists an edge $F$ of $\mathcal{H}$ that contains $\s{u,v,y}$.

Consider the edge $F'$ of $\mathcal{N}(\mathcal{H})$ that corresponds to the closed neighborhood of $y$ in the 2-section of $\mathcal{H}$. Then, $\s{u,v,y}\subseteq F'$. Given that $x\in I_{\mathcal{N}(\mathcal{H})}(uv)$, we have that $x\in F'$, which implies that $x$ is in the closed neighborhood of $y$ in the 2-section of $\mathcal{H}$, that is, there exists an edge of $\mathcal{H}$ that contains $\s{x,y}$.

Now suppose that $x\not\in I_{\mathcal{N}(\mathcal{H})}(uv)$. Then, there exists a vertex $y$ whose neighborhood in the 2-section of $\mathcal{H}$ has $u$ and $v$ but does not have $x$. Suppose without loss of generality that $v\in T[u,y]$. Consider an edge $F$ of $\mathcal{H}$ that contains $u$ and $y$. Since $T[F]$ is a subtree, we conclude that $v\in F$, so $\s{u,v,y}\subseteq F$. However, there is no edge of $\mathcal{H}$ containing $x$ and $y$.

Also note that the edges of $\mathcal{H}$ that contain $\s{u,v}$ are the same as the ones that contain $I_\mathcal{H}(uv)$, so the proof is complete.

\end{proof}

\begin{theorem}
  Let $\mathcal{H}$ be a hypertree and $T$ be a host tree of $\mathcal{H}$. Then, $\mathcal{H}$ is basic if and only if, for every edge $uv$ of $T$, the following holds: for every vertex $x$ not in $I_{\mathcal{H}}(uv)$, there exists a vertex $y$ such that there is an edge of $\mathcal{H}$ that contains $\s{u,v,y}$ and there is no edge of $\mathcal{H}$ that contains $\s{x,y}$.

  In other words, $\mathcal{H}$ is basic if and only if for every basic set $B$ of $\mathcal{H}$ and every vertex $x$ not in $B$, there exists a vertex $y$ such that there is an edge of $\mathcal{H}$ that contains $B\cup\s{y}$ and there is no edge of $\mathcal{H}$ that contains $\s{x,y}$.
\end{theorem}

\begin{proof}

We know from the previous lemma that every host tree of $\mathcal{H}$ is also a host tree of $\mathcal{N}(\mathcal{H})$. We infer from Lemma \ref{hosttcont} that $I_\mathcal{H}(uv)\subseteq I_{\mathcal{N}(\mathcal{H})}(uv)$ for every edge $uv$ of $T$.

On the other side, the condition that ``for every vertex $x$ not in $I_{\mathcal{H}}(uv)$, there exists a vertex $y$ such that there is an edge of $\mathcal{H}$ that contains $\s{u,v,y}$ and there is no edge of $\mathcal{H}$ that contains $\s{x,y}$'' is by Proposition \ref{basicn} equivalent to $I_{\mathcal{N}(\mathcal{H})}(uv)\subseteq I_\mathcal{H}(uv)$, which is in our context necessary and sufficient for the equality $I_\mathcal{H}(uv)=I_{\mathcal{N}(\mathcal{H})}(uv)$ to hold. Combine this with Theorem \ref{test} to conclude that the condition is necessary and sufficient for $\mathcal{H}$ to be a basic hypertree.

The second phrasing of the theorem in terms of basic sets follows immediately from the connection between the sets $I_\mathcal{H}(uv)$ and the basic sets of $\mathcal{H}$ given by Theorem \ref{wabasic}.

\end{proof}

It is interesting to note that, for two equivalent hypertrees, one could be basic unlike the other. Consider for example the hypertree $\mathcal{H}$ whose vertex set is $\s{1,2,3}$ and its edges are $\s{1,2}$ and $\s{2,3}$. Then, $\mathcal{N}(\mathcal{H})$ has the same edges, with the addition of the edge $\s{1,2,3}$, so $\mathcal{H}$ and $N(\mathcal{H})$ are equivalent and $\mathcal{H}$ is basic. However, $\s{1,2,3}$ is the only possible edge of $\mathcal{N}(\mathcal{N}(\mathcal{H}))$, so $\mathcal{N}(H)$ is not basic.

In \cite{correspondence}, \emph{basic chordal} graphs were defined as those chordal graphs whose clique trees are exactly the compatible trees of its clique graph. In other words, a chordal graph $G$ is basic chordal when $D\mathcal{C}(G)$ and $\mathcal{C}(K(G))$ are equivalent hypertrees. Much like there is a characterization of hypertrees in terms of the chordality of the line graph, there exists a similar characterization of basic hypertrees.

A hypergraph $\mathcal{H}$ is \emph{separating} if, for every two vertices $u$ and $v$ of it, there exists one edge that contains $u$ and does not contain $v$.

\begin{proposition} \cite{clique}\label{hellysep}
  Let $\mathcal{H}$ be a Helly and separating hypergraph. Then, the hypergraph of maximal cliques of $L(\mathcal{H})$ is equal of $D\mathcal{H}$.
\end{proposition}

\begin{proposition}
  Let $\mathcal{H}$ be a Helly and separating hypergraph. Then, its host trees are exactly the clique trees of $L(\mathcal{H})$.
\end{proposition}

\begin{proof}
  If $L(\mathcal{H})$ is not chordal, then $\mathcal{H}$ is not a hypertree, so neither $\mathcal{H}$ has a host tree nor $L(\mathcal{H})$ has a clique tree.

  Suppose now that $L(\mathcal{H})$ is chordal, so $\mathcal{H}$ is a hypertree by Theorem \ref{arbor}. The clique trees of $L(\mathcal{H})$ are just the host trees of the dual of the hypergraph of maximal cliques of $L(\mathcal{H})$, which by Proposition \ref{hellysep} is equal to the dual of $D\mathcal{H}$, and this is in turn isomorphic to $\mathcal{H}$. Therefore, the clique trees of $L(\mathcal{H})$ correspond to the host trees of $\mathcal{H}$.
\end{proof}

\begin{theorem}
  Let $\mathcal{H}$ be a hypergraph and $\mathcal{H}'$ be the hypergraph obtained from $\mathcal{H}$ by adding, for every $v\in V(\mathcal{H})$, the edge $\s{v}$, in case it is missing. Then, $\mathcal{H}$ is a basic hypertree if and only if it is Helly and  $L(\mathcal{H}')$ is basic chordal.
\end{theorem}

\begin{proof}

 Note that $\mathcal{H}'$ is a basic hypertree if and only if $\mathcal{H}$ is a basic hypertree, since the addition of unit sets is an operation that preserves the host trees and $\mathcal{N}(\mathcal{H})=\mathcal{N}(\mathcal{H}')$.

  In both directions of the proof we have that $\mathcal{H}$ (and hence $\mathcal{H}'$) is a hypertree, either from the definition of basic hypertree or, in the converse, if $\mathcal{H}$ is Helly and $L(\mathcal{H}')$ is basic chordal, then particularly $L(\mathcal{H}')$ is chordal and its induced subgraph $L(\mathcal{H})$ is also chordal, so we get that $\mathcal{H}$ is a hypertree from Theorem \ref{arbor}.


   Additionally, the presence of all the possible unit sets ensures that $\mathcal{H}'$ is separating, so the clique trees of $L(\mathcal{H}')$ correspond to the host trees of $\mathcal{H}'$. Given that the maximal cliques of $L(\mathcal{H}')$ are the edges of $D\mathcal{H}'$, it is not difficult to verify that the dually chordal graph that arises from taking the clique graph of $L(\mathcal{H}')$ is isomorphic to the 2-section of $\mathcal{H}'$, whose compatible trees are the trees on the same vertex set where every closed neighborhood induces a subtree, which are just the host trees of $\mathcal{N}(\mathcal{H}')$.

   As a consequence of this, if we are given that $\mathcal{H}$ is a basic hypertree, then $\mathcal{H}'$ is also a basic hypertree, $\mathcal{H}'$ and $\mathcal{N}(\mathcal{H}')$ have the same host trees and, from the previous paragraph, the clique trees of $L(\mathcal{H}')$ are the same as the compatible trees of $K(L(\mathcal{H}'))$, that is, $L(\mathcal{H}')$ is basic chordal.

   Conversely, if $L(\mathcal{H}')$ is basic chordal, then repeat the argument inversely to conclude that $\mathcal{H}$ is a basic hypertree.

\end{proof}

\section{The edges of host trees}

The main goal of this section is to show what edges can appear in a host tree of a hypertree and how these edges can be combined.

To determine the edges that can appear in a host tree, define for a hypergraph $\mathcal{H}$ and a set $A$ of vertices of $\mathcal{H}$ the hypergraphs $\mathcal{H}_A$ whose vertices are those of $\mathcal{H}$ and whose edges are the edges of $\mathcal{H}$ that contain $A$. On the other side, the \emph{complementary hypergraph} $\overline{\mathcal{H}_{A}}$ will also have the same vertices, but its edges are the edges of $\mathcal{H}$ that do not contain $A$. When $A=\s{u,v}$, the notation $\mathcal{H}_{uv}$ and $\overline{\mathcal{H}_{uv}}$ is used.


\begin{theorem} \label{isedge}
  Let $\mathcal{H}$ be a hypertree and $u$ and $v$ be two vertices of it. Then, $uv$ is the edge of some host tree of $\mathcal{H}$ if and only if $u$ and $v$ are in different connected components of the 2-section of $\overline{\mathcal{H}_{uv}}$.
\end{theorem}

\begin{proof}
  Suppose that $uv$ is the edge of a host tree $T$ of $\mathcal{H}$, which is also a host tree of $\overline{\mathcal{H}_{uv}}$. Let $T_u$ and $T_v$ be the two connected components of $T-uv$ containing $u$ and $v$, respectively. Since every edge of $\overline{\mathcal{H}_{uv}}$ is contained in $T_u$ or in $T_v$, the 2-section of $\overline{\mathcal{H}_{uv}}$ is disconnected, with $u$ and $v$ in different connected components of it.

  Conversely, suppose that $u$ and $v$ are in different connected components of the 2-section of $\overline{\mathcal{H}_{uv}}$ and let $T$ be a host tree of $\mathcal{H}$. Since the 2-section of $\overline{\mathcal{H}_{uv}}$ has no path from $u$ to $v$, $T[u,v]$ has an edge $e$ such that no edge of $\overline{\mathcal{H}_{uv}}$ contains $e$. Thus, every edge of $\mathcal{H}$ that contains $e$ also contains $uv$. By Lemma \ref{-+}, $T-e+uv$ is a host tree of $\mathcal{H}$ that contains the edge $uv$.
\end{proof}

With Theorem \ref{isedge}, some results about the edges of clique trees of chordal graphs and the compatible trees of dually chordal graphs become corollaries.

\begin{corollary} \cite{correspondence,reduced}
  Let $G$ be a chordal graph, and let $C$ and $C'$ be two maximal cliques of $G$. Let $S=C\cap C'$. Then, there exists a clique tree of $G$ that contains the edge $CC'$ if and only if $G-S$ is not connected, with the vertices of $C\setminus C'$ and $C'\setminus C$ in different connected components. In that case, we additionally have that, given the hypertree $\mathcal{H}$ that is the dual of the hypergraph of maximal cliques of $G$, the basic set $I_\mathcal{H}(CC')$ consists of the maximal cliques of $G$ that contain $S$.
\end{corollary}

\begin{proof}
  Every clique tree of $G$ is a host tree of $\mathcal{H}$. Note that $\mathcal{H}_{CC'}$ has the edges of the form $\mathcal{C}_v$, where $v\in S$, so the edges of $\overline{\mathcal{H}_{CC'}}$ are of the same form but with $v\not\in S$. On the other side, it is not difficult to verify that $C$ and $C'$ are not adjacent in the 2-section of $\overline{\mathcal{H}_{CC'}}$ and every other pair $C_1,C_2$ of maximal cliques of $G$ is adjacent in the 2-section of $\overline{\mathcal{H}_{CC'}}$ if and only if $C_1\cap C_2\cap (V(G)\setminus S)\neq\emptyset$.

  Let $u\in C\setminus C'$ and $w\in C'\setminus C$. If there exists an $uw$-path in $G-S$ with edges $e_1,...,e_k$, let $C_1,...,C_k$ be maximal cliques of $G$ containing each of those edges. Thus, $C$ and $C'$ can be connected in the 2-section of $\overline{\mathcal{H}_{CC'}}$ through a walk that uses these cliques. Conversely, if there is a $CC'$-path in the 2-section of $\overline{\mathcal{H}_{CC'}}$, then taking for every pair of consecutive cliques of that path a vertex in their intersection that is not in $S$ gives a sequence of vertices that can be used to connect $u$ and $w$ in $G-S$.

  Combine the previous two paragraphs with Theorem \ref{isedge} to obtain the first part of the corollary.

  To verify that $I_\mathcal{H}(CC')$ consists of the maximal cliques of $G$ that contain $S$, we have just to apply the definition of this set together with the form we found for the edges of $\mathcal{H}_{CC'}$.
\end{proof}

\begin{corollary}\label{isedgedu}
  Let $G$ be a dually chordal graph and $u$ and $v$ be two vertices of $G$. Let $G'$ be the subgraph of $G$ (with the same vertex set) such that every edge $xy$ of $G$ is also in $G'$ if and only if $\s{x,y}\neq\s{u,v}$ and there exists $z\not\in N[u]\cap N[v]$ such that $\s{x,y}\subseteq N[z]$. Then, $uv$ is the edge of some compatible tree of $G$ if and only if there exists an $uv$-path in $G'$. In that case, we additionally have that, given the hypertree $\mathcal{H}$ of maximal cliques of $G$, the basic set $I_\mathcal{H}(uv)$ consists of the vertices $w$ such that $N[u]\cap N[v]\subseteq N[w]$.
\end{corollary}

\begin{proof}
  Every compatible tree of $G$ is a host tree of the hypertree $\mathcal{H}$ of maximal cliques of $G$. Additionally, $\overline{\mathcal{H}_{uv}}$ consists of all the maximal cliques of $G$ that do not contain $\s{u,v}$.

  Let us verify that the 2-section of $\overline{\mathcal{H}_{uv}}$ is just equal to the subgraph $G'$. Let $C$ be a maximal clique of $G$ that does not contain $\s{u,v}$ and let $x$ and $y$ be two elements of it. Since we do not have both $u$ and $v$ in $C$, there exists a vertex $z$ in $C$ that is not adjacent to both $u$ and $v$. This vertex $z$ also satisfies that $\s{x,y}\subseteq N[z]$. Thus, every edge of the 2-section of $\overline{\mathcal{H}_{uv}}$ is an edge of $G'$.

  Conversely, let $xy$ be an edge of $G'$ and let $z\not\in N[u]\cap N[v]$ be such that $\s{x,y}\subseteq N[z]$. Let $C$ be a maximal clique of $G$ that contains $\s{x,y,z}$. By the definition of $z$, $C$ does not contain $\s{u,v}$, so $C$ is an edge of $\overline{\mathcal{H}_{uv}}$ and $xy$ is an edge of its 2-section.

  We can finish the proof of the first part of the corollary applying Theorem \ref{isedge}.

  To verify that $I_\mathcal{H}(uv)$ consists of the vertices $w$ such that $N[u]\cap N[v]\subseteq N[w]$, note that $\mathcal{H}$ is equivalent to the hypertree $\mathcal{H}'$ of closed neighborhoods of $G$, so $I_\mathcal{H}(uv)=I_{\mathcal{H'}}(uv)$. By definition, $I_{\mathcal{H'}}(uv)$ is the intersection of all the closed neighborhoods of vertices in $N[u]\cap N[v]$, which consists just of all the vertices $w$ such that $N[u]\cap N[v]\subseteq N[w]$.

\end{proof}

The first part of Corollary \ref{isedgedu} is a result that has not been published before, while the second part appears in \cite{correspondence}.

The previous results gave conditions about the structure of the hypertree or of the graph, depending on the case, but we can also establish conditions in terms of host trees, in case one is known.

\begin{theorem} \label{givenT}
  Let $\mathcal{H}$ be a hypertree and $T$ be a host tree of it. Let $u$ and $v$ be two different vertices of $\mathcal{H}$. Then, $uv$ is the edge of some host tree of $\mathcal{H}$ if and only if there exists an edge $e$ in $T[u,v]$ such that $I_\mathcal{H}(uv)=I_\mathcal{H}(e)$.
\end{theorem}

\begin{proof}
  Suppose that there exists an edge $e$ in $T[u,v]$ such that $I_\mathcal{H}(uv)=I_\mathcal{H}(e)$. By Lemma \ref{-+}, $T-e+uv$ is a host tree of $\mathcal{H}$.

  Conversely, suppose that $uv$ is the edge of some host tree of $\mathcal{H}$. If for every edge $e$ of $T[u,v]$ there exists an edge of $\mathcal{H}$ that contains $e$ and does not contain $\s{u,v}$, a contradiction to Theorem \ref{isedge} arises, as $T[u,v]$ would be an $uv$-path in the 2-section of $\overline{\mathcal{H}_{uv}}$. Thus, there exists an edge $e$ of $T[u,v]$ such that every edge of $\mathcal{H}$ that contains $e$ also contains $\s{u,v}$. On the other side, as $T$ is a host tree of $\mathcal{H}$, every edge of $\mathcal{H}$ that contains $\s{u,v}$ also contains $e$. Therefore, $I_\mathcal{H}(uv)=I_\mathcal{H}(e)$.
\end{proof}

Note however that the condition that $I_\mathcal{H}(uv)$ is a basic set of $\mathcal{H}$ is not sufficient so that $\mathcal{H}$ has a host tree with the edge $uv$. Consider for example the hypertree $\mathcal{H}$ with vertices 1, 2, 3 and 4 and with edges $\s{1,2}$, $\s{2,3}$ and $\s{1,2,3,4}$, whose basic sets are just the same sets. We have that $I_\mathcal{H}(13)=I_\mathcal{H}(14)=\s{1,2,3,4}$. Every host tree of $\mathcal{H}$ has the edges 12 and 23 and one edge between 4 and another vertex, so 14 is the edge of some host tree but 13 is not.

As before, Theorem \ref{givenT} can be rephrased to obtain the particular results for chordal and dually chordal graphs.

\begin{corollary}
\hfill
  \begin{enumerate}
    \item \cite{reduced} Let $G$ be a chordal graph, $T$ be a clique tree of $G$ and $C_1$ and $C_2$ be two different maximal cliques of $G$. Then, $C_1C_2$ is the edge of some clique tree of $G$ if and only if there exists an edge $C_3C_4$ in $T[C_1C_2]$ such that $C_1\cap C_2=C_3\cap C_4$.
    \item Let $G$ be a dually chordal graph, $T$ be a compatible tree of $G$ and $u$ and $v$ be two different vertices of $G$. Then, $uv$ is the edge of some compatible tree of $G$ if and only if there exists an edge $xy$ in $T[u,v]$ such that $N[u]\cap N[v]=N[x]\cap N[y]$.
  \end{enumerate}
\end{corollary}

\begin{proof}
  Apply Theorem \ref{givenT} to the dual of the hypergraph of maximal cliques of $G$ to obtain 1. Apply Theorem \ref{givenT} to the hypergraph of closed neighborhoods of the vertices of $G$ to obtain 2.
\end{proof}

Suppose now that we have two host tree edges $uv$ and $xy$ such that both $I_\mathcal{H}(uv)$ and $I_H(xy)$ are equal to the same basic set $B$. As a consequence of Proposition \ref{every}, the edges of $\mathcal{H}$ that contain $\s{u,v}$ are the same as the edges that contain $\s{x,y}$, which are just all the edges that contain $B$. For that reason, $\mathcal{H}_{uv}=\mathcal{H}_{xy}=\mathcal{H}_B$. The same equality holds for their complementary hypergraphs and their 2-sections. Keeping this in mind, we can reformulate Theorem \ref{isedge}.

\begin{lemma}\label{bsubihe}
  Let $\mathcal{H}$ be a hypegraph, and $A$ and $B$ be two subsets of $V(\mathcal{H})$. If $B$ is not contained in a single connected component of $2S(\overline{\mathcal{H}_A})$, then $A\subseteq I_\mathcal{H}(B)$.
\end{lemma}

\begin{proof}
  If there exists an edge $E$ of $\mathcal{H}$ that contains $B$ but does not contain $A$, then $E$ ensures that $B$ is a clique of $2S(\overline{\mathcal{H}_A}))$, which is not possible by our hypothesis. Thus, every edge of $\mathcal{H}$ that contains $B$ also contains $A$ and hence $A\subseteq I_\mathcal{H}(B)$.
\end{proof}

\begin{theorem}\label{lasdeb}
  Let $\mathcal{H}$ be a hypertree and $B$ be a basic set of $\mathcal{H}$. Let $u$ and $v$ be two different elements of $B$. Then, $I_{\mathcal{H}}(uv)=B$ and $uv$ is the edge of some host tree of $\mathcal{H}$ if and only if $u$ and $v$ are in two different connected components of the 2-section of $\overline{\mathcal{H}_B}$.
\end{theorem}

\begin{proof}
Suppose first that $I_{\mathcal{H}}(uv)=B$ and $uv$ is the edge of some host tree of $\mathcal{H}$. That $I_{\mathcal{H}}(uv)=B$ implies that $\overline{\mathcal{H}_{uv}}=\overline{\mathcal{H}_B}$ so, by Theorem \ref{isedge}, that $uv$ is the edge of some host tree implies that $u$ and $v$ are in two different connected components of  $2S(\overline{\mathcal{H}_B})$.

Conversely, suppose that $u$ and $v$ are in two different connected components of the 2-section of $\overline{\mathcal{H}_B}$. By Lemma \ref{bsubihe}, $B\subseteq I_\mathcal{H}(uv)$, while it is trivial that every edge of $\mathcal{H}$ that contains $B$ also contains $u$ and $v$, so $I_\mathcal{H}(uv)\subseteq I_\mathcal{H}(B)$ by Proposition \ref{every}. Since $B=I_\mathcal{H}(xy)$ for some edge $xy$, we have the equality $I_\mathcal{H}(B)=B$. Therefore, $I_\mathcal{H}(uv)=B$.

We also infer from the previous paragraph that $\overline{\mathcal{H}_B}$ and $\overline{\mathcal{H}_{uv}}$ are equal, and so are their 2-sections. Thus, $u$ and $v$ are in two different connected components of $2S(\overline{\mathcal{H}_{uv}})$ and, by Theorem \ref{isedge}, $uv$ is the edge of a host tree $T$ of $\mathcal{H}$.
\end{proof}



We will see now that similar ideas can also be applied to obtain a new characterization of hypertrees.

\begin{lemma}\label{component}
  Let $\mathcal{H}$ be a hypergraph and $A$ be a set of vertices of $\mathcal{H}$. Then, every connected component of $\overline{\mathcal{H}_{A}}$ is in $Comp(\mathcal{H})$.
\end{lemma}

\begin{proof}
  Every connected component of $\overline{\mathcal{H}_{A}}$ has a single vertex or is the connected union of the edges of $\mathcal{H}$ that are contained in it, so it is in $Comp(\mathcal{H})$.
\end{proof}

\begin{theorem}\label{newchyp}
  Let $\mathcal{H}$ be a hypergraph. Then, $\mathcal{H}$ is a hypertree if and only if, for every basic set $B$ of $\mathcal{H}$, $B$ induces a disconnected graph in $2S(\overline{\mathcal{H}_B})$.
\end{theorem}

\begin{proof}
  Suppose that $\mathcal{H}$ is a hypertree and let $B$ be a basic set of $\mathcal{H}$. We know that $B=I_\mathcal{H}(uv)$ for some edge $uv$ of a host tree $T$ of $\mathcal{H}$. By Theorem \ref{lasdeb}, $u$ and $v$ are vertices of $B$ in different connected components of $2S(\overline{\mathcal{H}_B})$.

  The proof of the converse will be by induction on the amount of vertices that $\mathcal{H}$ has.

  If $\mathcal{H}$ has no more than two vertices, then it is trivial that $\mathcal{H}$ is a hypertree. On the other side, $\mathcal{B}(\mathcal{H})$ is either empty or the only basic set of $\mathcal{H}$ is $V(\mathcal{H})$, which satisfies the condition of the theorem.

  Suppose now that the implication is true for every hypergraph with at most $k$ vertices and let $\mathcal{H}$ have $k+1$ vertices.
  
  Let $B$ be any basic set of $\mathcal{H}$. Then (reasoning like in the proof of Lemma \ref{bsubihe}) every edge of $\mathcal{H}$ contains $B$ or is contained in a connected component of $2S(\overline{\mathcal{H}_B})$. It is easy to verify that the nonempty intersection (or connected union) of two sets such that each either contains $B$ or is contained in a connected component of  $2S(\overline{\mathcal{H}_B})$ also satisfies that it contains $B$ or is contained in a connected component of $2S(\overline{\mathcal{H}_B})$. Therefore, every edge of $Comp(\mathcal{H})$ contains $B$ or is contained in a connected component of $2S(\overline{\mathcal{H}_B})$, so $Comp(\mathcal{H})$ satisfies the hypothesis of the theorem.
  
  Fix now $B$ as a maximal basic set of $\mathcal{H}$ and let $B'$ be any other basic set of $\mathcal{H}$. Since it is not a possibility that $B\subseteq B'$, reasoning like in the previous paragraph we derive that $B'$ is contained in some connected component of $2S(\overline{\mathcal{H}_B})$.
  
  Let $A$ be any connected component of $2S(\overline{\mathcal{H}_B})$. If $A$ has a single vertex $v$, then let $T_A$ be the tree whose only vertex is $v$.
  
  Consider now the case that $A$ has two or more vertices. Let $\mathcal{C}^A$ be the hypergraph with vertex set $A$ and whose edges are those of $Comp(\mathcal{H})_A$. Note that $Comp(\mathcal{C}^A)=\mathcal{C}^A$, so the basis of $\mathcal{C}^A$ consists of the basic sets of $\mathcal{H}$ that are contained in $A$.
  
  Since $Comp(\mathcal{H})$ satisfies the hypothesis of the theorem, so does $\mathcal{C}^A$, which has fewer basic sets and edges. Thus, by the induction hypothesis, $\mathcal{C}^A$ is a hypertree. Let $T_A$ be a host tree of $\mathcal{C}^A$.
  
  Combine all the trees $T_A$ and, in case $i$ of these trees ($i\geq 2$) have vertices of $B$, add $i-1$ edges between vertices of $B$ to ensure that all the trees $T_A$ such that $A\cap B\neq\emptyset$ form a tree. Finally, add edges arbitrarily to connect with the trees $T_A$, where $A\cap B=\emptyset$, until the result is a tree $T$.

  Let us see that every basic set $B'$ of $\mathcal{H}$ induces a subtree of $T$. If $B'\neq B$, then there exists a connected component $A$ of $2S(\overline{\mathcal{H}_B})$ that contains $B'$, so $B'$ induces a subtree of $T_A$ and hence of $T$.

  Suppose now that $B'=B$. For every connected component $A$ of $2S(\overline{\mathcal{H}_B})$ such that $A\cap B\neq\emptyset$, the set $A\cap B$ is by Lemma \ref{component} an edge of $Comp(\mathcal{H})$, and hence of $\mathcal{C}^A$, so $T[A\cap B]$ is equal to $T_A[A\cap B]$ and it is a subtree. The edges that we added during the construction of $T$ ensure that $T[B]$ is a subtree.

  Therefore, every basic set of $\mathcal{H}$ induces a subtree of $T$, so this is a host tree of $\mathcal{H}$ and $\mathcal{H}$ is a hypertree.

\end{proof}

Given a basic set $B$ of a hypertree $\mathcal{H}$, define $A(B)$ as the set whose elements are the connected components of the 2-section of $\overline{\mathcal{H}_{B}}$ that have at least one vertex of $B$ and denote its cardinality by $\alpha_{B}$.

\begin{theorem} \label{numberedgesT}
  Let $\mathcal{H}$ be a hypertree, $T$ be a host tree of $\mathcal{H}$ and $B$ be a basic set of $\mathcal{H}$. Then, the number of edges $e$ in $T$ such that $I_\mathcal{H}(e)=B$ is equal to $\alpha_{B}-1$.
\end{theorem}

\begin{proof}

  Suppose that $B$ has nonempty intersection with $k$ connected components $A_1,...,A_k$ of the 2-section of $\overline{\mathcal{H}_{B}}$. For every $i$ between 1 and $k$, the set $A_i\cap B$ is an edge of $Comp(\mathcal{H})$ by Lemma \ref{component}. The subtree $T[B]$ is formed by the disjoint subtrees $T[A_i\cap B]$, $1\leq i \leq k$, and $k-1$ more edges $e_1,...,e_{k-1}$ connecting them. For every $i$ between 1 and $k-1$, the endpoints of $e_i$ are elements of $B$, so by part 1 of Proposition \ref{completion} we have that $I_\mathcal{H}(e_i)\subseteq B$. On the other side, the inclusion $B\subseteq I_\mathcal{H}(e_i)$ is inferred from Lemma \ref{bsubihe}. We conclude that $B= I_\mathcal{H}(e_i)$.

  Consider now any edge $e$ of $T$ such that $I_\mathcal{H}(e)=B$. If it is not contained in $B$, then the equality $B=I_\mathcal{H}(e)$ cannot hold as $I_\mathcal{H}(e)$ does contain $e$. Thus, $e$ must be contained in $B$. By Theorem \ref{lasdeb}, $e$ can only be one of $e_1,...,e_{k-1}$.

  
  Therefore, the set of edges $e$ of $T$ such that $I_\mathcal{H}(e)= B$ consists just of the edges $e_1,...,e_{k-1}$.
\end{proof}

We also state the following property that is inferred from the first paragraph of the proof of Theorem \ref{numberedgesT}.

\begin{proposition} \label{gamma}
  Let $\mathcal{H}$ be a hypertree, $T$ be a host tree of $\mathcal{H}$ and $B$ be a basic set of $\mathcal{H}$. Consider the graph $\Gamma_{B,T}$ whose vertex set is $A(B)$ and such that two different connected components $A_1$ and $A_2$ in $A(B)$ are adjacent if and only if $T$ has an edge $e$ such that $I_\mathcal{H}(e)=B$ and $e$ has one endpoint in $A_1$ and the other endpoint in $A_2$. Then, $\Gamma_{B,T}$ is a tree.
\end{proposition}

Given a hypertree $\mathcal{H}$ and a basic set $B$ of it, let $\Delta(B)$ consist of the edges $uv$ contained in $B$ such that $u$ and $v$ are in different connected components of the 2-section of $\overline{\mathcal{H}_B}$. By Theorem \ref{lasdeb}, every edge $e$ in $\Delta(B)$ is the edge of some host tree of $\mathcal{H}$ and $I_\mathcal{H}(e)=B$; actually, $\Delta(B)$ contains all the edges satisfying these two conditions. Let us now say that a collection $E$ of edges in $\Delta(B)$ is $B$-\emph{admissible} if the graph $\Gamma_{B,E}$ that has vertex set $A(B)$ and for every $e$ in $E$ there is an edge in $\Gamma_{B,E}$ between the connected components in $A(B)$ that contain each endpoint of $e$ is a tree.

The following is a result that generalizes a characterization of clique trees that appears in \cite{tesis}.

\begin{theorem} \label{admissible}
  Let $\mathcal{H}$ be a hypertree whose basic sets are $B_1,..,B_m$ and let $T$ be a tree with vertex set $V(\mathcal{H})$. Then, $T$ is a host tree of $\mathcal{H}$ if and only if $E(T)=\cup_{i=1}^m{E_i}$, where $E_i$ is $B_i$-admissible for $1\leq i\leq m$.
\end{theorem}

\begin{proof}
  Assume that $T$ is a host tree of $H$. That $E(T)$ is the union of admissible sets is a direct consequence of Proposition \ref{gamma}.

  Conversely, assume now that $E(T)=\cup_{i=1}^m{E_i}$ where, for $1\leq i\leq m$, $E_i$ is $B_i$-admissible. Let $T'$ be a host tree of $\mathcal{H}$ such that $E(T')=\cup_{i=1}^m{E_i'}$, where every $E_i'$ is $B_i$-admissible. Given any $i$ between 1 and $m$, let us see that the tree $T''$ such that $T''=T'-E_i'+E_i$ is a host tree of $H$. Suppose to the contrary that $T''$ is not a host tree of $\mathcal{H}$. Since $\mathcal{H}$ is equivalent to $\mathcal{B}(Comp(\mathcal{H}))$, there is a basic set $B$ such that $T''[B]$ is not a subtree and, under this condition, suppose that $B$ is inclusionwise minimal.

  For every connected component $A_1$ in $A(B)$ we have that $A_1\cap B$ is in $Comp(\mathcal{H})$, as it is the intersection of two sets in $Comp(\mathcal{H})$. If $A_1\cap B$ has just one vertex, then it clearly induces a subtree in $T'$. If it has more vertices, then it is the connected union of one or more basic sets that are strictly contained in $B$ and, by the minimality condition of $B$, these basic sets induce a subtree in $T''$ and hence $T''[A_1\cap B]$ is a subtree.

  Let $B=B_j$, for some $j$ between 1 and $m$. Combine the previous paragraph with the fact that $E_j$ and $E_j'$ are $B_j$-admissible to conclude that $T''[B]$ is a subtree, thus contradicting what we had initially assumed.

  Therefore, $T''[B]$ is a subtree for every basic set of $\mathcal{H}$ and $T''$ is a host tree of $\mathcal{H}$. We can apply this fact, starting with the host tree $T'$, to successively replace every $E_i'$ with $E_i$ to conclude that $T$ is also a host tree of $\mathcal{H}$.

\end{proof}

We will now use this characterization to show that edges associated to the same basic set are interchangeable to obtain one host tree from another.

\begin{proposition}
  Let $\mathcal{H}$ be a hypertree, $B$ be a basic set of $\mathcal{H}$ and $e$ and $e'$ be two edges whose endpoints are in different connected components in $A(B)$. Then, there exists a host tree $T$ that has the edge $e$ and such that $T-e+e'$ is also a host tree.
\end{proposition}

\begin{proof}
  Let $A_1$ and $A_2$ be the different connected components in $A(B)$ that contain the endpoints of $e$ and $A_3$ and $A_4$ be the different connected components in $A(B)$ that contain the endpoints of $e'$. Suppose without loss of generality that, in case of equality, the possible equalities between these components are $A_1=A_3$ and $A_2=A_4$. Consider a total order of the elements of $A(B)$ where $A_3\leq A_1 < A_2\leq A_4$ and $A_1$ and $A_2$ are consecutive in the order. For every two consecutive elements of $A(B)$ in the total order, pick one element from each of them to form an edge, with the restriction that for the case of $A_1$ and $A_2$ that edge has to be $e$. Then, the collection $E$ of those edges is $B$-admissible and so is $(E\setminus\s{e})\cup\s{e'}$. Let $T$ be a host tree of $\mathcal{H}$ whose $B$-admissible set is $E$. Then, $T-e+e'$ will have $B$-admissible set $(E\setminus\s{e})\cup\s{e'}$ while all the other admissible sets remain the same. Therefore, $T-e+e'$ is a host tree of $\mathcal{H}$.
\end{proof}

Denote the collection of all the host trees of a hypergraph $\mathcal{H}$ by $\tau(\mathcal{H})$. We know how to express a basic set of a hypertree as an intersection of edges. Now we will see a way to express it as a union.

\begin{proposition}\label{swap}
  Let $\mathcal{H}$ be a hypertree and $B$ be a basic set of $\mathcal{H}$. Let $u$ and $v$ be elements of $B$ that are in different connected components of $2S(\overline{\mathcal{H}_B})$. Then, $B=\bigcup_{T\in\tau(\mathcal{H})}{T[u,v]}$.
\end{proposition}

\begin{proof}
  Since $B$ induces a subtree of every host tree of $\mathcal{H}$, we have that $T[u,v]\subseteq B$ for every host tree $T$ of $\mathcal{H}$. Thus, $\bigcup_{T\in\tau(\mathcal{H})}{T[u,v]}\subseteq B$.

  Suppose now that $B$ has an element $w$ different from $u$ and $v$, and also suppose without loss of generality that $w$ is not in the same connected component of $2S(\overline{\mathcal{H}_B})$ as $u$. Then, by Proposition \ref{swap}, there exist host trees $T_1$ and  $T_2$ of $\mathcal{H}$ such that $T_2=T_1-uv+uw$, so $w\in T_2[u,v]$. We conclude from this argument that $B\subseteq\bigcup_{T\in\tau(\mathcal{H})}{T[u,v]}$.
\end{proof}

\section{Host trees as maximum weight spanning trees}

Let $G$ be a graph. We say that $G$ is a \emph{weighted} graph if it comes with a function $w$ that assigns to each edge of $G$ a real number, which is called the weight of the edge. Suppose that $G$ is connected. Given a spanning tree $T$ of $G$, the weight of $T$ is defined as $w(T)=\sum_{e\in E(T)}{w(e)}$. The problem of the \emph{minimum weight spanning tree} consists in finding a spanning tree of $G$ with minimum weight. There are several algorithms to find a minimum weight spanning tree, the most known ones being Prim's algorithm \cite{prim} and Kruskal's Algorithm \cite{kruskal}. The problem of the maximum weight spanning tree can be defined analogously and can be solved by making minor adjustments to the cited algorithms.

The goal of this short section is to revisit some properties of maximum weight spanning trees to find the connection between them and host trees of hypergraphs.

A weighted graph may or may not have a unique maximum weight spanning tree. In the following, we will consider the case that there is more than one maximum weight spanning tree.

\begin{proposition}\label{swap2}
  Let $T$ and $T'$ be two different maximum weight spanning trees of a weighted graph $G$. Then, there exists an edge $e\in E(T)\setminus E(T')$ and an edge $e'\in E(T')\setminus E(T)$ such that $T-e+e'$ is also a maximum weight spanning tree of $G$.
\end{proposition}

\begin{proof}
  Among all the edges in the symmetric difference $E(T)\Delta E(T')$, let $e$ be of minimum weight and suppose without loss of generality that $e\in E(T)$. Since $T'$ is connected, it contains an edge $e'$ whose endpoints are in different connected components of $T-e$. We cannot have the inequality $w(e')<w(e)$ by the way $e$ was defined. The inequality $w(e)<w(e')$ is not true either, otherwise the spanning tree $T-e+e'$ would have larger weight than $T$. Thus, $w(e)=w(e')$, and the tree $T-e+e'$ satisfies that $w(T-e+e')=w(T)$, so $T-e+e'$ is also a maximum weight spanning tree of $G$.
\end{proof}

\begin{corollary}\label{sequence}
Let $T$ and $T'$ be two different maximum weight spanning trees of a weighted graph $G$ such that $|E(T)\setminus E(T')|=k$. Then, there exists a sequence $T_1T_2....T_{k+1}$ of maximum weight spanning trees of $G$ such that $T_1=T$, $T_{k+1}=T'$ and, for $1\leq i\leq k$, $T_{i+1}=T_i-e_i+e'_i$, where $e_i\in E(T)$ and $e'_i\in E(T')$. Furthermore, there is no shorter sequence satisfying the same conditions.
\end{corollary}

\begin{proof}
  Apply Proposition \ref{swap2} successively starting with $T$ until all edges of $E(T)\setminus E(T')$ are removed and all edges of $E(T')\setminus E(T)$ are added to get $T'$ at the end. This cannot be done in fewer steps as the $k$ edges of $E(T')\setminus E(T)$ need to be added and we cannot add more than one edge per step.
\end{proof}

\begin{corollary} \label{multiplicity}
  Let $T$ and $T'$ be two maximum weight spanning trees of a weighted graph $G$. Then, the sets of edge weights of $T$ and $T'$ are the same. What is more, each weight has the same multiplicity in $T$ and $T'$.
\end{corollary}

\begin{proof}
The procedure of Corollary \ref{sequence} ensures that in each step the weights and their multiplicities stay the same.
\end{proof}

This result resembles Theorem \ref{numberedgesT}. Actually, we will see that is not by chance, as host trees are particular cases of maximum weight spanning trees.

\begin{lemma}\label{wineq}
  Let $\mathcal{H}$ be a hypergraph. Consider the complete graph with vertex set $V(\mathcal{H})$ and, for every edge $uv$ of it, we give it the weight $w(uv)=|\{F\in E(\mathcal{H}):~\{u,v\}\subseteq F\}|$ and let $T$ be a spanning tree. Then, $w(T)\leq\sum_{F\in E(\mathcal{H})}{|F|} - |E(\mathcal{H})|$.
\end{lemma}

\begin{proof}
\hfill
  \[\sum_{e\in E(T)}{w(e)}=\sum_{e\in E(T)}|\s{F\in E(\mathcal{H}):~e\subseteq F}|=\sum_{e\in E(T)}\sum_{F\in E(\mathcal{H})}|\s{F}\cap\s{F'\in E(\mathcal{H}):~e\subseteq F'}|=\]
  \[\sum_{F\in E(\mathcal{H})}\sum_{e\in E(T)}|\s{F}\cap\s{F'\in E(\mathcal{H}):~e\subseteq F'}|=\sum_{F\in E(\mathcal{H})}|E(T[F])|\leq \sum_{F\in E(\mathcal{H})}(|F|-1)=\]
  \[\sum_{F\in E(\mathcal{H})}{|F|} - |E(\mathcal{H})|\]
\end{proof}

\begin{theorem} \label{hostmax}
   Let $\mathcal{H}$ be a hypergraph. Consider the complete graph with vertex set $V(\mathcal{H})$ and, for every edge $uv$ of it, we give it the weight $w(uv)=|\{F\in E(\mathcal{H}):~\{u,v\}\subseteq F\}|$. Then, $\mathcal{H}$ is a hypertree if and only if  every maximum weight spanning tree has weight $\sum_{F\in E(\mathcal{H})}{|F|} - |E(\mathcal{H})|$. What is more, the host trees of $\mathcal{H}$ are all the spanning trees of weight $\sum_{F\in E(\mathcal{H})}{|F|} - |E(\mathcal{H})|$.
\end{theorem}

\begin{proof}
A spanning tree $T$ is a host tree of $\mathcal{H}$ if and only if $T[F]$ is a subtree for every edge $F$ of $\mathcal{H}$, that is, if and only if $|E(T[F])|=|F|-1$ for every edge $F$ of $\mathcal{H}$. Thus, the inequality of the proof of Lemma \ref{wineq} becomes an equality and a spanning tree is a host tree if and only if its weight is $\sum_{F\in E(\mathcal{H})}{|F|} - |E(\mathcal{H})|$, and hence a maximum weight spanning tree. This proves the second statement of the theorem, from which the first one follows.

\end{proof}

Given that the host trees of a hypertree can be presented as maximum weight spanning trees, the question arises whether the maximum weight spanning trees of every complete graph correspond to the host trees of some hypergraph.

Consider the complete graph with four vertices 1, 2, 3 and 4 such that $w(12)=w(23)=w(34)=w(14)=2$ and $w(13)=w(24)=1$. In this case, the maximum weight spanning trees are the spanning trees of the cycle $12341$. Suppose that these four spanning trees are all the host trees of a hypertree $\mathcal{H}$. By proposition \ref{swap}, the set \s{1,2,3,4} would be the only basic set of $\mathcal{H}$, so every spanning tree would be a host tree, a contradiction. Therefore, not every set of maximum weight spanning trees of a complete weighted graph is the set of host trees of some hypertree.

Like in previous sections, there are some particular results about chordal and dually chordal graphs that can now be presented as special cases of Theorem \ref{hostmax}.

\begin{theorem} \cite{mckee}
  Let $G$ be a chordal graph. Consider the complete graph whose vertices are the maximal cliques of $G$, and give each edge $CC'$ of it weight $|C\cap C'|$. Then, the clique trees of $G$ are the maximum weight spanning trees of this graph. Additionally, the weight of every clique tree is equal to $\sum_{v\in V(G)}{|\mathcal{C}_v|}~-~|V(G)|$.
\end{theorem}

\begin{proof}
  Let $\mathcal{H}$ be the dual of the hypergraph of maximal Cliques of $G$. The clique trees of $G$ are the host trees of $\mathcal{H}$, and the given weighting matches the one induced by $\mathcal{H}$, so Theorem \ref{hostmax} can be applied to finish the proof.
\end{proof}

Particularly, if $G$ is connected, it is not difficult to verify that the clique trees of $G$ are spanning trees of $K(G)$, so we can refine the result by considering the maximum weight spanning trees of $K(G)$ under the same weighting.

\begin{theorem}
  Let $G$ be a dually chordal graph and consider the complete graph on the same vertices as $G$. Consider the following two weighings for this graph.

  \begin{enumerate}
    \item Every edge $uv$ has the amount of maximum cliques of $G$ that contain $\s{u,v}$ as weight.
    \item \cite{dra} Every edge $uv$ has $|N_G[u]\cap N_G[v]|$ as weight.
  \end{enumerate}

Then, the compatible trees of $G$ are the maximum weight spanning trees of $G$ under any of the two weightings. In the first case, the total weight is $\sum_{C\in \mathcal{C}(G)}|C|~-~|\mathcal{C}(G)|$, while in the second case the total weight is $2|E(G)|$.
\end{theorem}

\begin{proof}
  Let $\mathcal{H}_1$ be the hypergraph of maximal cliques of $G$, and $\mathcal{H}_2$ be the hypergraph of closed neighborhoods of the vertices of $G$. Apply Theorem \ref{hostmax} to $\mathcal{H}_1$ to conclude that the compatible trees are the maximum weight spanning trees with respect to the first weighting. Apply Theorem \ref{hostmax} to $\mathcal{H}_2$ to conclude that the compatible trees are the maximum weight spanning trees with respect to the second weighting.
\end{proof}

Particularly, if $G$ is connected, one can verify that the compatible trees of $G$ are spanning trees of $G$, so it is possible to restrict the weighting to $G$ to find all the compatible trees.

We end the section with an analog of the initial results.

\begin{proposition}
Let $T$ and $T'$ be two different host trees of a hypertree $\mathcal{H}$. Then...

\begin{enumerate}
  \item There exist edges $e\in E(T)\setminus E(T')$ and $e'\in E(T')\setminus E(T)$ such that $T-e+e'$ is also a host tree.
  \item If $|E(T)\setminus E(T')|=k$, then there exists a sequence $T_1T_2....T_{k+1}$ of host trees of $\mathcal{H}$ such that $T_1=T$, $T_{k+1}=T'$ and, for $1\leq i\leq k$, $T_{i+1}=T_i-e_i+e'_i$, where $e_i\in E(T)$ and $e'_i\in E(T')$. Furthermore, there is no shorter sequence satisfying the same conditions.
\end{enumerate}
\end{proposition}

\begin{proof}
  By Theorem \ref{hostmax}, there exists a weighted complete graph whose maximum weight spanning trees are the host trees of $\mathcal{H}$. Then, everything follows from Proposition \ref{swap2} and Corollary \ref{sequence}.
\end{proof}

\section{Conclusions}

To date, most of the knowledge about the structure of host trees of hypertrees came from studying their most notable exponents, namely, the clique trees of chordal graphs and the compatible trees of dually chordal graphs.

The main goal of this paper was introducing a general theory about the host trees of hypertrees, with proofs as short as possible that only require basic hypergraph notions. A second goal was to give the reader an opportunity to read about the host trees of hypertrees in a single place, with a unified notation and without the need to read several papers that use different approaches and terminology.

Theorems \ref{wabasic}, \ref{test}, \ref{isedge}, \ref{givenT}, \ref{lasdeb}, \ref{numberedgesT}, \ref{admissible} and \ref{hostmax} had similar versions that were exclusively for clique trees of chordal graphs, which were transformed here into hypertree results whose proofs have weak or null links with the theory of chordality.

Furthermore, novel elements like the concept of equivalent hypergraphs and the new characterization of hypertrees in Theorem \ref{newchyp} appear. The concept of basic hypertree, although inspired in basic chordal graphs, is also completely new.

For future work, it would be interesting to consider other known graph representations and study them from a hypergraph perspective.


\begin{thebibliography}{99}

\bibitem{fagin2} C. Beeri, R. Fagin, D. Maier, M. Yannakakis, On the desirability of acyclic database schemes, Journal ACM 30 (1983) 479–513.

\bibitem{algor} A. Brandstäd, V. Chepoi, F. Dragan, The algorithmic use of hypertree structure and maximum neighbourhood orderings, Discrete Applied Mathematics 82 (1998), 43–77.
    
\bibitem{brands} A. Brandst\"{a}dt, F. Dragan, V. Chepoi and V. Voloshin, Dually chordal graphs, SIAM J. Discrete Math.  11  (1998), 437–455.

\bibitem{modeling} Dai, Q., Gao, Y. (2023). Hypergraph Modeling. In: Hypergraph Computation. Artificial Intelligence: Foundations, Theory, and Algorithms. Springer, Singapore.

\bibitem{tesis} P. De Caria, A joint study of chordal and dually chordal graphs, Ph. D. Thesis, Universidad Nacional de La Plata, 2012.

\href{https://www.mate.unlp.edu.ar/~pdecaria/decariathesis.pdf}{https://www.mate.unlp.edu.ar/~pdecaria/decariathesis.pdf}

\bibitem{braz} P. De Caria, M. Gutierrez, Determining what sets of trees can be the clique trees of a chordal
graph, Journal of the Brazilian Computer Society 18 (2012) 121–128.


\bibitem{correspondence} P. De Caria, M. Gutierrez, On the correspondence between tree representations of chordal and
dually chordal graphs, Discrete Applied Mathematics 164 (2014) 500–511.

\bibitem{dra} F. Dragan, HT-graphs: centers, connected r-domination and Steiner trees, Comput. Sci. J. Moldova 1 (1993) 64–83.

\bibitem{fagin} R. Fagin, Acyclic database schemes of various degrees: a painless introduction, Lecture Notes in Comput. Sci. 159 (1983), 65–89.

\bibitem{fagin3} R. Fagin, Degrees of acyclicity for hypergraphs and relational database schemes, J. ACM 30 (1983) 514–550.

\bibitem{arbor} C. Flament, Hypergraphes arbores, Discrete Mathematics, 21 (1978), 223–227.

\bibitem{gavril} F. Gavril, The intersections of subtrees in trees are exactly the chordal graphs, J. Combin. Theory Ser. B 116 (1974) 47–56.

\bibitem{clique} M. Gutierrez, J. Meidanis, Algebraic theory for the clique operator, J. Braz. Comp. Soc. 7 (2001) 53–64.

\bibitem{reduced} M. Habib, J. Stacho, Reduced clique graphs of chordal graphs, European Journal of Combinatorics 33 (2012) 712–735.

\bibitem{kruskal} J.B. Kruskal, On the shortest spanning tree of a graph and the traveling salesman problem, Proceedings of the American Mathematical Society 7 (1956) 48–50.

\bibitem{retrieval} Q. Liu, Y. Huang,  D.N. Metaxas, Hypergraph with sampling for image retrieva", Pattern Recognition 44 (2011) 2255–2262.

\bibitem{mckee} T.A. McKee, How chordal graphs work, Bulletin of the ICA 9 (1993) 27–39.

\bibitem{bioinf} R. Patro, C. Kingsoford, Predicting protein interactions via parsimonious network history inference, Bioinformatics, 29 (2013) 237–246.

\bibitem{prim} R. C. Prim, Shortest Connection Networks and some Generalizations, Bell System Technical Journal 36 (1957) 1389–1401.


\bibitem{spin} J. P. Spinrad, Efficient graph representations, Field Institute Monographs, American Mathematical Society, Providence, Rhode Island, 2003.

\bibitem{enum} P. S. Kumara, C.E. Veni Madhavanb, Clique tree generalization and new subclasses of chordal graphs, Discrete Applied Mathematics 117 (2002) 109–131.


\end{thebibliography}
\end{document}